\documentstyle[11pt,amssymb]{amsart}

\newcommand{\seq}[3]{{#1}_{#2}, \ldots, {#1}_{#3}}
\newcommand{\Dws}{D^{(w,\Sigma)}}
\newcommand{\Y}{\Sigma \times {{\Bbb S}}^1}
\newcommand{\Spz}{\text{Sp}\, (2g,{\Bbb Z})}
\newcommand{\Ct}{{\Bbb C}[[t]]}
\newcommand{\ima}{{\bf i}}
\newcommand{\whff}{\widetilde{HFF}{}^*_g}
\newcommand{\rhff}{\overline{HFF}{}^*_g}

\newcommand{\la}{\langle}
\newcommand{\ra}{\rangle}

\newcommand{\inc}{\hookrightarrow}
\newcommand{\ar}{\rightarrow}
\newcommand{\bd}{\partial}
\newcommand{\x}{\times}
\newcommand{\ox}{\otimes}
\newcommand{\iso}{\cong}

\newcommand{\point}{\text{pt}}
\newcommand{\CP}{{\Bbb C \Bbb P}}

\newcommand{\Sym}{\text{Sym}}

\newcommand{\PD}{\text{P.D.}}

\renewcommand{\AA}{{\Bbb A}}
\newcommand{\CC}{{\Bbb C}}
\newcommand{\DD}{{\Bbb D}}

\newcommand{\QQ}{{\Bbb Q}}

\newcommand{\SS}{{\Bbb S}}
\newcommand{\ZZ}{{\Bbb Z}}

\renewcommand{\a}{\alpha}
\renewcommand{\b}{\beta}
\renewcommand{\d}{\delta}
\newcommand{\g}{\gamma}

\renewcommand{\l}{\lambda}

\newcommand{\p}{\phi}
\newcommand{\q}{\psi}

\renewcommand{\S}{\Sigma}
\newcommand{\D}{\Delta}
\renewcommand{\L}{\Lambda}

\theoremstyle{plain}
\begingroup
\theorembodyfont{\sl}
\newtheorem{thm}{Theorem}
\newtheorem{cor}[thm]{Corollary}
\newtheorem{lem}[thm]{Lemma}
\newtheorem{prop}[thm]{Proposition}

\endgroup

\theoremstyle{definition}
\newtheorem{defn}[thm]{Definition}

\theoremstyle{remark}
\newtheorem{rem}[thm]{Remark}

\title{Basic classes for four-manifolds not of simple type}

\author{Vicente Mu\~noz}
\address{Departamento de \'Algebra Geometr\'{\i}a y Topolog\'{\i}a \\
 Facultad de Ciencias \\ Universidad de M\'alaga \\ 29071 M\'alaga \\ Spain}
\email{vmunoz@@agt.cie.uma.es}
\thanks{\hbox{$^*$}Partially supported by DGES through Spanish Research Project
 PB97-1095. 
 }
\date{November, 1998}

\begin{document}

\maketitle

\begin{abstract}
  We extend the notion of basic classes (for the Donaldson invariants) to
  $4$-manifolds with $b^+>1$ which are (potentially) not of simple type
  or satisfy $b_1 >0$. We also give
  a structure theorem for the Donaldson invariants of $4$-manifolds
  with $b^+>1$, $b_1 > 0$ and of strong simple type.
\end{abstract}

\section{Introduction}
\label{sec:intro}

The primary purpose of this paper is to extend the notion of basic classes (as
defined by Kronheimer and Mrowka~\cite{KM} for the
Donaldson invariants) to all $4$-manifolds satisfying $b^+>1$. So we drop the
conditions $b_1=0$ and being of simple type imposed in~\cite{KM}. We give a
(somewhat partial) structure theorem for the Donaldson invariants of any
$4$-manifold with $b^+>1$ (compare~\cite{Kr}).

Donaldson invariants for a (smooth, compact, oriented) $4$-manifold $X$
with $b^+>1$ (and with a homology orientation) are defined as linear
functionals~\cite{KM}
$$
  D^w_X: \AA(X)= \Sym^*(H_0(X) \oplus H_2(X)) \ox \L^* H_1(X) \ar \CC,
$$
where $w \in H^2(X;\ZZ)$. As for the grading of
$\AA(X)$, we give degree $4-i$ to the elements in $H_i(X)$.
If $x \in H_0(X)$ denotes the class of a point,
we say that $X$ is of $w$-simple type when $D^w_X((x^2-4)z)=0$,
for any $z\in \AA(X)$. When $X$ has $b_1=0$ and it
is of $w$-simple type for some $w \in H^2(X;\ZZ)$, then it is of
$w'$-simple type for any other $w'\in H^2(X;\ZZ)$, and it is
said to be of simple type for brevity~\cite{FS}.
For $4$-manifolds of simple type (with $b^+ >1$),
it is customary to define~\cite{KM} the formal
power series (in a variable $t$) given by
$$
\DD_X^{w}(tD) =D_X^w((1+{x \over 2}) e^{tD}),
$$
for any $D \in H_2(X)$. We introduce the following notions as in~\cite{hff}

\begin{defn}
\label{def:type}
  Let $w \in H^2(X;\ZZ)$. We say that $X$ is of {\bf $w$-finite type} when
there
  is some $n \geq 0$ such that $D^w_X((x^2-4)^nz)=0$, for any $z\in \AA(X)$.
  The order (of $w$-finite type) is the minimum of such $n$.
  $X$ is of finite type if it is of $w$-finite type for any $w\in H^2(X;\ZZ)$
and the order
  (of finite type) is the supremum of all the orders of $w$-finite type of $X$,
  with $w\in H^2(X;\ZZ)$.

  \noindent Also $X$ is of {\bf $w$-strong simple type}, for $w \in
H^2(X;\ZZ)$, when
  $D^w_X((x^2-4)z)=0$, for any $z\in \AA(X)$ and $D^w_X(\d z)=0$, for any $\d
\in H_1(X)$ and
  any $z\in \AA(X)$.
  $X$ is of strong simple type if it is of $w$-strong simple type for any $w\in
H^2(X;\ZZ)$.
\end{defn}

Clearly strong simple type implies simple type, and the two conditions are
equivalent
when $b_1=0$. Also finite type of order $1$ means simple type and order $0$
means that
the Donaldson invariants are identically zero. Note that by~\cite[theorem
7.3]{hff},
all $4$-manifolds with $b^+>1$ are of finite type.
Since $D_X^{w+2\a}=(-1)^{\a^2}D_X^w$, for any $\a\in H^2(X;\ZZ)$, the order
of $w$-finite type only depends on $w$ modulo $2H^2(X;\ZZ)$, and hence the
order of finite type is finite. We shall see later
(theorem~\ref{thm:orderft}) that the order of $w$-finite type does not depend
on
$w\in H^2(X;\ZZ)$ and hence the order of finite type is finite and equal to any
order of $w$-finite type.

Our first result is an extension of the structure theorem for the
Donaldson invariants as given in~\cite{KM}~\cite{FS} to $4$-manifolds of
strong simple type, in the form that was
stated in~\cite[proposition 7.6]{hff}
(recall that $Q$ stands for the intersection form).

\begin{thm}
\label{thm:sst}
  Let $X$ be a manifold with $b^+ >1$ and of $w$-strong simple type
  for some $w \in H^2(X;\ZZ)$. Then $X$ is of strong simple type and
  $\DD_X^{w}= e^{Q /2} \sum (-1)^{{K_i \cdot {w} +{w}^2} \over 2} a_i \,
  e^{K_i}$, for finitely many $K_i \in H^2(X;\ZZ)$ and rational numbers
  $a_i$ (the collection is empty
  when the invariants all vanish). These classes are lifts
  to integral cohomology of ${w}_2(X)$. Moreover, for any embedded
  surface $S \inc X$ of genus $g$, representing a non-torsion
  homology class and with $S^2 \geq 0$,
  one has $2g-2 \geq S^2 +|K_i \cdot S|$, for all $K_i$.
\end{thm}

When we put $b_1=0$ we recover the structure theorem of
Kronheimer and Mrowka given in~\cite{KM} for the Donaldson invariants
of $4$-manifolds of simple type with $b_1=0$ and $b^+ >1$ (see also~\cite{FS}
for the case of simply connected $4$-manifolds).
It is possible that the proof provided in~\cite{KM} can be adapted to prove
theorem~\ref{thm:sst}, but there are difficulties as one of the starting
conditions
for the analysis in~\cite{KM} is that of admissibility
in~\cite[definition 2.23]{KM}, which can not be removed in general.
A.\ Stipsicz informed the author that he encountered this same problem in the
computation of the Donaldson invariants of the $4$-torus~\cite{Stipsicz}.
We have decided to include a self-contained proof of theorem~\ref{thm:sst}
as it is used by the author in the determination of the Donaldson
invariants of the product of two compact surfaces of genus bigger than one
in~\cite[subsection 7.3]{hff}. It will also serve to clarify the ideas
behind the proof of our more general structure theorem for
the Donaldson invariants of arbitrary $4$-manifolds with
$b^+>1$ (theorem~\ref{thm:basic}).

\begin{defn}
\label{def:basic-sst}
  If $X$ is a $4$-manifold with $b^+ >1$ and of strong simple type, then the
  {\bf basic classes} of $X$ will be the cohomology classes $K_i$ of
  theorem~\ref{thm:sst} such that $a_i \neq 0$.
\end{defn}

\begin{rem}
\label{rem:sst}
  There are many $4$-manifolds with $b_1 >0$ and $b^+>1$ which are of
  strong simple type. For instance, $X=\S_g \x \S_h$, where $\S_r$ is a
  compact surface of genus $r$ and $g, h\geq 1$, are of strong simple
  type, as it is proved in~\cite[theorem 7.11]{hff}. On the other hand, there
  are
  also $4$-manifolds with $b^+>1$ not of strong simple type. If $X$ is
  any $4$-manifold with $b^+>1$ and non-vanishing Donaldson invariants,
  then $X' =X \# \, \SS^1 \x \SS^3$ is not of strong simple type, as
  follows from~\cite[lemma 5.8]{hff}.
  Nonetheless for $X$ of strong simple type and $X' =X \#\, \SS^1 \x \SS^3$,
  let $\d$ be the image in $H_1(X')$ of the generator of
  $\pi_1(\SS^1 \x \SS^3)$. We have
  $$
  D^w_{X'}(\d \, (1+{x \over 2}) e^{tD})= D_X^w ((1+{x \over 2}) e^{tD})
  =\DD^w_X(tD),
  $$
  which has the same shape as that provided by theorem~\ref{thm:sst}.
\end{rem}

For general $4$-manifolds with $b^+>1$ it is not licit to suppose that
$(x^2-4)$ and all $\d\in H_1(X)$ kill the Donaldson invariants. In principle,
there might be $4$-manifolds
with $b_1=0$, $b^+>1$ not of simple type~\cite{Kr}
(although they have not been found so far) and, in any case,
remark~\ref{rem:sst}
tells us that there are examples where the $1$-homology classes act
non-trivially.
We have the following two main results:

\begin{thm}
\label{thm:orderft}
  Let $X$ be a $4$-manifold with $b^+>1$. Let $w,w'\in H^2(X;\ZZ)$. Then the
  order of $w$-finite type and the order of $w'$-finite type of $X$ are equal.
  This number is thus the order of finite type of $X$.
\end{thm}

\begin{thm}
\label{thm:basic}
  Let $X$ be a $4$-manifold with $b^+>1$, $w \in H^2(X;\ZZ)$ and $z \in
\L^*H_1(X) \subset
  \AA(X)$ an homogeneous element. Then there are
  finitely many cohomology classes $K_i \in H^2(X;\ZZ)$ and non-zero
  polynomials $p_i, q_i \in \Sym^* H^2(X) \ox \QQ[\l]$ such that
  \begin{equation}
    D^w_X(z e^{tD+\l x})= e^{Q(tD)/2+2\l}\sum p_i(tD,\l) e^{K_i\cdot tD} +
     e^{-Q(tD)/2-2\l}\sum q_i(tD,\l) e^{\ima K_i\cdot tD},
  \label{eqn:falta}
  \end{equation}
  for any $D\in H_2(X)$. The collection of classes $K_i$ is independent of
  $w$ and $K_i$ are lifts to integral cohomology of ${w}_2(X)$.
\end{thm}

\begin{rem}
\label{rem:symm}
  Let $d_0=-w^2-{3\over 2}(1-b_1+b^+)$ and $d=\deg (z)/2$. Here
  $d_0-d \in \ZZ$. Then
  $D_X^w(z  e^{\ima tD - \l x})= \ima^{d_0-d} D_X^w(z e^{tD+\l x})$,
  $\ima=\sqrt{-1}$, so
  $q_i(tD,\l)=\ima^{d-d_0}
  p_i(\ima tD, -\l)$ for all $i$. Therefore the $p_i$ determine the
  $q_i$ and conversely.
  Also if $K_j=-K_i$ then $p_j(-tD,\l)= (-1)^{d_0-d}p_i(tD,\l)$. So the classes
  $K_i$ come actually in pairs $\pm K_i$.
\end{rem}

We finally have the following definition, which agrees with
definition~\ref{def:basic-sst}
for strong simple type manifolds.

\begin{defn}
\label{def:basic}
  Let $X$ be a $4$-manifold with $b^+ >1$ and $z \in \L^*H_1(X)$ an
  homogeneous element. The cohomology classes $K_i$ of theorem~\ref{thm:basic}
  such that $p_i \neq 0$ are called {\bf basic classes} for $(X,z)$.
  The union of all the basic classes for $(X,z)$ where $z$ runs through
  any homogeneous basis of $\L^*H_1(X)$ is the set of basic classes for $X$.
\end{defn}

The proof of theorem~\ref{thm:basic} is based on two techniques. On the one
hand, the
Fukaya-Floer homology of the three-manifold $Y=\Y$, where $\S$ is a
compact surface, as determined in~\cite{hff}, which we recall in
section~\ref{sec:hff}
for the convenience of the reader. On the other hand, partial
use of the blow-up formula in~\cite{FS-bl} which relates the Donaldson
invariants of a $4$-manifold $X$ and those of its blow-up
$\tilde X=X \# \overline{\CP}^2$ .
Section~\ref{sec:sst} is devoted to prove theorem~\ref{thm:sst} and in
section~\ref{sec:basic} we study the Donaldson invariants of general
$4$-manifolds with $b^+>1$ in order to prove theorems~\ref{thm:orderft}
and~\ref{thm:basic}.

\section{Fukaya-Floer homology revisited}
\label{sec:hff}

Let $Y$ be a $3$-manifold with $b_1 >0$ and $w \in H^2(Y;\ZZ/2\ZZ)$ non-zero.
For
any loop $\d \subset Y$, we have defined the Fukaya-Floer homology~\cite{hff}
$$
  HFF_*(Y,\d),
$$
which is a $\Ct$-module, endowed with a $\Ct$-bilinear pairing
$$
\la\/,\ra:HFF_*(Y,\d) \otimes HFF_*(- Y, -\d) \ar \Ct,
$$
where $-Y$ is $Y$ with reversed orientation. For every $4$-manifold $X_1$ with
boundary $\bd X_1=Y$ and $w_1 \in H^2(X;\ZZ)$ such that $w_1|_Y=w \in
H^2(Y;\ZZ/2\ZZ)$, $z \in \AA(X_1)$ and $D_1 \subset X_1$ a $2$-cycle
with $\bd D_1=\d$, one has a relative invariant
$$
  \p^w(X_1,e^{tD_1}) \in HFF_*(Y,\d).
$$
The relevant gluing theorem is:

\begin{thm}[\cite{hff}]
\label{thm:gluing}
  Let $X=X_1 \cup_Y X_2$ and $w \in H^2(X;\ZZ)$ such that there exists
  $\S \in H^2(X;\ZZ)$ whose Poincar\'e dual lies in the image of
  $H_2(Y;\ZZ) \ar H_2(X;\ZZ)$,
  and satisfies $w \cdot \S \equiv 1\pmod 2$. Put
  $w_i =w|_{X_i} \in H^2(X_i;\ZZ)$.
  Let $D \in H_2(X)$ decomposed as $D=D_1 +D_2$ with $D_i \subset
  X_i$, $i=1,2$, $2$-cycles with $\bd D_1=\d$, $\bd D_2=-\d$. Choose
  $z_i \in \AA(X_i)$, $i=1,2$. Then
  $$
     \Dws_X(z_1z_2e^{tD})=
     \la\p^{w_1}(X_1,z_1e^{tD_1}),\p^{w_2}(X_2,z_2e^{tD_2})\ra,
  $$
  where $\Dws_X=D^w_X+D^{w+\S}_X$.
  When $b^+=1$ the invariants are calculated for metrics on $X$ giving a
  long neck.
\end{thm}

We restrict to the case $Y=\Y$, where $\S$ is a surface of genus $g\geq 1$,
$w=\PD [\SS^1]$, $\d=\SS^1 \subset Y$. As $Y \iso (-Y)$, we have a natural
identification $HFF_*(Y,\SS^1) \iso HFF^*(Y,\SS^1)$ and a pairing
$$
   \la\/,\ra:HFF^*(Y,\SS^1) \otimes HFF^*( Y, \SS^1) \ar \Ct.
$$
Let $A=\S \x D^2$ be the product of $\S$ times a $2$-dimensional disc and
consider the horizontal section $\D=\S\x\point \subset A$.
Let $\{\g_i\}$ be a symplectic basis for $H_1(\S)$
with $\g_i \cdot \g_{g+i}=1$, $1\leq i \leq g$.
The Fukaya-Floer homology $HFF^*(Y,\SS^1)$
is actually a $\Ct$-algebra~\cite[section 5]{hff} generated by
$\a=2\p^w(A,\S e^{t\D})$, $\b=-4 \p^w(A,x e^{t\D})$ and
$\q_i=\p^w(A,\g_i e^{t\D})$, $1\leq i \leq 2g$, where the product
is determined by $\p^w(A, z_1e^{t\D})\p^w(A, z_2 e^{t\D})=
\p^w(A, z_1z_2 e^{t\D})$, $z_1, z_2\in \AA(\S)$. In particular
$\p^w(A, z e^{t\D})$, $z \in\AA(\S)$, generate $HFF^*(Y,\SS^1)$.
The mapping class group of $\S$ acts on $HFF^*(Y,\SS^1)$ factoring through an
action of the symplectic group $\Spz$ on $\{\q_i\}$.
The invariant part is generated by $\a$, $\b$
and $\g= -2 \sum \q_i\q_{g+i}$.

  In general, we shall have the following situation: $X$ is a $4$-manifold with
  $b^+ >1$ and $\S \inc X$ is an embedded surface of genus $g\geq 1$ and
  self-intersection
  zero such that there exists $w \in H^2(X;\ZZ)$ with $w\cdot \S \equiv 1\pmod
  2$. Let $A$ be a small tubular neighbourhood of $\S$ in $X$ and
  denote by $X_1$ the complement of the interior of $A$, so that $X_1$ is a
  $4$-manifold with boundary $\bd X_1=Y=\Y$ and $X=X_1 \cup_Y A$.
  Take $D \in H_2(X)$ with $D\cdot \S =1$. Then $D$ can be represented by a
  cycle $D=D_1 + \D$, where $D_1 \subset X_1$, $\bd D_1 =\SS^1$ and
  $\D=\point\x D^2 \subset A$. Theorem~\ref{thm:gluing} says
  $$
  \Dws_X(z e^{tD +s\S+\l x})=\la\p^w(X_1,z e^{tD_1+s\S+\l x}),\p^w(A,e^{t\D})\ra,
  $$
  for any $z \in\AA(X_1)$. Here $t$ is a formal variable, and both terms are
  meant to be developed in powers of $s$ and $\l$. Also $\S$ denotes
  the homology class represented by the surface and the corresponding
  cohomology class by Poincar\'e duality.
  The effective Fukaya-Floer homology, $\whff$, (see~\cite[definition 5.10]{hff})
  is the sub-$\Ct$-module of $HFF^*_g$ with the property that
  $$
  \p^w(X_1,z e^{tD_1}) \in \whff \subset HFF^*_g.
  $$
  If $R$ is a polynomial such that $R(\a,\b)=0$ acting on $\whff$ then
  $$
  R(2{\bd\over \bd s}, -4{\bd\over \bd \l}) \p^w(X_1,ze^{tD_1+s\S+\l x})=0,
  $$
  and therefore $R(2{\bd\over \bd s}, -4{\bd\over \bd \l}) \Dws_X(z e^{tD
  +s\S+\l x})=0$.
  The Donaldson invariants $\Dws_X$ satisfy thus partial differential equations
  coming from the polynomials vanishing on $\whff$. We have the following result

\begin{prop} {\rm (\cite[theorem 5.13]{hff})}
\label{prop:HFF}
  The effective Fukaya-Floer homology is a sub-$\Ct$-module $\whff \subset
  HFF^*_g$ such that there is a direct sum decomposition
  $$
   \whff= \bigoplus_{r=-(g-1)}^{g-1} R_{g,r},
  $$
  where $R_{g,r}$ are free $\Ct$-modules such that, in $R_{g,r}$,
  $\a- (4r\ima -2t)$, $\b-8$ and
  $\g$ are nilpotent if $r$ is even, $\a-( 4r+2t)$, $\b+8$ and $\g$
  are nilpotent if $r$ is odd.
\end{prop}

\section{Strong simple type manifolds}
\label{sec:sst}

This section is devoted to the proof of theorem~\ref{thm:sst}.
We start with two technical lemmas.
We shall abbreviate strong simple type to ``sst'' in this section.

\begin{lem}
\label{lem:1}
  Let $X$ be a $4$-manifold with $b^+>1$ and  $w,\S \in H^2(X;\ZZ)$
  such that $\S^2=0$ and $w \cdot \S \equiv 1 \pmod 2$.
  If $X$ is both of $w$-sst and of $(w+\S)$-sst and $D\in H_2(X)$
  then there exist
  power series $f_{r,D}(t)$, $-(g-1) \leq r \leq g-1$, such that
\begin{equation}
\label{eqn:DD}
  \DD^w_X (tD +s\S)=e^{Q(tD +s\S)/2}\sum\limits_{r=-(g-1)}^{g-1}
  f_{r,D}(t) e^{2rs},
\end{equation}
  i.e. $\DD^w_X (tD +s\S)$ is a solution of the differential equation
  $\prod_{r=-(g-1)}^{g-1} ({\bd \over \bd s}-(2r +t(D\cdot \S)))$. Moreover
  $$
  \DD^{w+\S}_X (tD +s\S)=e^{Q(tD +s\S)/2}\sum\limits_{r=-(g-1)}^{g-1}
  (-1)^{r+1} f_{r,D}(t) e^{2rs}.
  $$
\end{lem}

\begin{pf}
  It is enough to prove~\eqref{eqn:DD} for $D \in H_2(X)$ with $D\cdot\S=1$,
  using linearity and continuity.
  Consider an embedded surface $\S\inc X$ of genus $g \geq 1$ representing
  the Poincar\'e dual of $\S \in H^2(X;\ZZ)$. Then we are in the situation
  described in section~\ref{sec:hff} and will stick to the notations used
  therein. The relative Donaldson invariant
  $$
    \p^w(X_1,e^{tD_1}) \in HFF_*(\Y,\SS^1)
  $$
  lies in the kernels of $\b^2-64$ and all $\q_i$, $1 \leq i \leq 2g$. This is
  clear
  as, for instance, $\la(\b^2-64) \p^w(X_1,e^{tD_1}),\p^w(A, z e^{t\D})\ra=
  \Dws_X(16 (x^2-4) z e^{tD})=0$, for any $\p^w(A, z e^{t\D})$,
  since $X$ is both of $w$-sst and of $(w+\S)$-sst.
  Thus $(\b^2-64) \p^w(X_1,e^{tD_1})=0$.

  Thus in order to compute
  $\Dws_X(e^{tD+s\S+\l x})=\la \p^w(X_1,e^{tD_1}),\p^w(A, e^{t\D+s\S+\l
  x})\ra$,
  it only matters the projection of $\p^w(A, e^{t\D+s\S+\l x})$ to
  the reduced Fukaya-Floer homology (see~\cite[definition 5.6]{hff})
  $$
   \rhff=HFF^*_g/(\b^2-64, \seq{\q}{1}{2g})HFF^*_g.
  $$
  By~\cite[theorem 5.9]{hff},
$$
  \rhff= \bigoplus\limits_{r=-(g-1)}^{g-1} \bar{R}_{g,r}, \qquad \text{where }
  \bar{R}_{g,r}=
  \left\{ \begin{array}{ll} \Ct[\a,\b] / (\a- (4r\ima-2t) ,\b-8) & \text{$r$
      even} \\  \Ct[\a,\b] / (\a- (4r+2t),\b+8) \qquad & \text{$r$ odd} 
  \end{array} \right.
$$
  Now we translate the relations of $\rhff$ into partial differential equations
  satisfied by the Donaldson invariants. For instance, the relation
  $$
  (\b-8) \prod_{-(g-1)\leq r\leq g-1 \atop \text{$r$ odd}} (\a-(4r+2t))=0
  $$
  gives the differential equation
  $$
  ({\bd \over\bd \l}+2) \prod_{-(g-1)\leq r\leq g-1 \atop \text{$r$ odd}} ({\bd
\over
   \bd s} -(2r+t)) \Dws_X(e^{tD +s\S+\l x})=0.
  $$
  This finally yields the existence of
  power series $g_r(t)$, $-(g-1) \leq r \leq g-1$ (we drop the $D$ from the
subindex),
  such that
  $$
   \Dws_X(e^{tD +s\S+\l x})=  \sum_{\text{$r$ odd}} g_r(t) e^{st +2rs +2\l}+
    \sum_{\text{$r$ even}} g_r(t) e^{-st +2r\ima s -2\l}.
  $$
  If $d_0=d_0(w) = -w^2 - {3\over 2}(1-b_1+b^+)$ denotes
  half the dimensions (modulo $4$) of the moduli
  spaces of anti-self-dual connections of $SO(3)$-bundles determined by $w$,
then
  $$
   D^w_X(e^{tD +s\S+\l x})= {1\over 2} \big( \sum_{\text{$r$ odd}}
   g_r(t) e^{st +2rs +2\l}+ \sum_{\text{$r$ even}} g_r(t) e^{-st +2r\ima s
-2\l} +
  $$
  $$
   + \sum_{\text{$r$ odd}} \ima^{-d_0}g_r(\ima t) e^{-st +2r\ima s -2\l}+
    \sum_{\text{$r$ even}} \ima^{-d_0}g_r(\ima t) e^{st -2rs +2\l} \big),
  $$
   and hence
  $$
   \DD^w_X(tD +s\S)=  \sum f_r(t) e^{st +2rs}, \qquad \text{where } f_r(t)=
   \left\{ \begin{array}{ll} g_r(t) & \text{$r$ odd} \\
   \ima^{-d_0}g_{-r}(\ima t) \qquad &
   \text{$r$ even} \end{array} \right.
  $$
  We leave $\DD^{w+\S}_X$ to the reader upon noting that $d_0(w+\S) \equiv
d_0(w)+2 \pmod 4$.
\end{pf}

\begin{lem}
\label{lem:2}
  In the situation of lemma~\ref{lem:1}, $X$ is of
  $w$-sst if and only if it is of $(w+\S)$-sst.
\end{lem}

\begin{pf}
  Arguing by contradiction, suppose $X$ is of $w$-sst but not of $(w+\S)$-sst.
  Then there exists $z=(x^2-4)^r \d_1 \cdots \d_p$ with $\d_i \in H_1(X)$ and
  $r+p >0$ such that $D^w_X(z e^{tD +s\S+\l x})=0$,
  $D^{w+\S}_X(z e^{tD +s\S+\l x})$ is non-zero but
  $D^{w+\S}_X(z (x^2-4) e^{tD +s\S+\l x})=0$ and
  $D^{w+\S}_X(z \d e^{tD +s\S+\l x})=0$ for any $\d \in H_1(\S)$.
  Keeping the notations of the proof of the previous lemma,
  $\p^w(X_1,z e^{tD_1})$
  lies in the kernel of $\b^2-64$ and $\g$ and so, arguing
  as in lemma~\ref{lem:1}, it is
  $$
   \Dws_X(z e^{tD +s\S+\l x})=  \sum_{\text{$r$ odd}} g_r(t) e^{st +2rs +2\l}+
    \sum_{\text{$r$ even}} g_r(t) e^{-st +2r\ima s -2\l},
  $$
  for some power series $g_r(t)$, $-(g-1) \leq r \leq g-1$.
  In particular $D^w_X(z e^{tD +s\S+\l x})$ is non-zero, so
  $X$ is not of $w$-sst. This contradiction proves the lemma.
\end{pf}

\noindent {\em Proof of theorem~\ref{thm:sst}.\/}
Now we proceed to the proof of theorem~\ref{thm:sst} by steps.

\noindent {\underline{Step 1}}. $X$ is of strong simple type.

Let $S$ be a $4$-manifold with $b^+>1$ and $w \in H^2(S;\ZZ)$. Let $\tilde S=
S\# \overline{\CP}^2$ denote its blow-up and let $E$ stand for the cohomology
class of the exceptional divisor. Therefore $H^2(\tilde S ;\ZZ)=
H^2(S;\ZZ) \oplus \ZZ E$.
The general blow-up formula~\cite{FS-bl} implies that $S$ is
of $w$-sst if and only if $\tilde S$ is of $w$-sst
if and only if $\tilde S$ is of $(w+E)$-sst.

With this said, suppose $X$ is of $w$-sst for some $w \in H^2(X;\ZZ)$.
We shall prove that $X$ is of $w'$-sst for any other $w' \in H^2(X;\ZZ)$.
Consider any cohomology class $L \in H^2(X;\ZZ)$ with $N=L^2>0$. Blow up
$X$ at $N$ points to obtain $\tilde X= X \# N \overline{\CP}^2$, with
$\seq{E}{1}{N}$ denoting the exceptional divisors. Let $L'=L -E_1
- \ldots -E_N \in H^2(\tilde X ;\ZZ)$ which has $(L')^2=0$. As $X$ is
of $w$-sst, $\tilde X$ will be both of $w$-sst and $(w+E_1)$-sst (recall
that $N >0$). Now $w \cdot L'$ and $(w+E_1)\cdot L'$ have different
parity since $E_1\cdot L'=1$. Therefore one of them is odd. If $w \cdot
L' \equiv 1 \pmod 2$, then lemma~\ref{lem:2} implies that $\tilde X$ is of
$(w+L')=(w+L-E_1 -\ldots -E_N)$-sst, and hence $X$ is of $(w+L)$-sst.
Alternatively, if $(w+E_1) \cdot L' \equiv 1 \pmod 2$, then $\tilde X$ is of
$(w+E_1+L')=(w+L-E_2 -\ldots -E_N)$-sst and $X$ of $(w+L)$-sst again.
In conclusion, $X$ is of $(w+L)$-sst, for any $L \in H^2(X;\ZZ)$ with $L^2>0$.

Now given $w$ and $w'$, take $T \in H^2(X;\ZZ)$ with $T^2>0$. For $n$ large,
it will be $(w'-w +n T)^2>0$. Considering $L=w'-w+n T$, it follows
that $X$ is of $(w+L)=(w'+nT)$-sst. Now taking $L=-n T$, we have that $X$
is of $w'$-sst, as required.

\noindent {\underline{Step 2}}.
There exists at least one $w\in H^2(X;\ZZ)$ with
$\DD_X^{w}= e^{Q /2} \sum a_{i,w} \, e^{K_i}$,
for finitely many cohomology classes $K_i \in
H^2(X;\ZZ)$ and rational numbers $a_{i,w}$.

Let $S$ be a $4$-manifold with $b^+>1$ and of sst. If $\tilde S$ is
the blow-up of $S$ with exceptional divisor $E$, then the blow-up
formula~\cite{FS-bl} says that $\DD_S^w(tD)=\DD_{\tilde S}^w (tD)$
and $\DD_S^{w}(tD)={\bd\over \bd r}|_{r=0} \DD_{\tilde S}^{w+E} (tD+rE)$, for
$D \in H_2(X)$, so we see that it is enough to prove the claim for $\tilde S$.

After possibly blowing up, we can suppose that $X$ has an
indefinite odd intersection form of the form $Q=r (1) \oplus s(-1)$,
with $r, s \geq 2$, $n=r+s$. Put $\seq{A}{1}{r},\seq{B}{1}{s}$
for the corresponding basis. Then we set
$\S_1=A_2-B_1$,
$\S_j=A_j +B_1$, $2 \leq j \leq r$,
$\S_{r + 1}=B_2-A_1$,
$\S_{r + j}=B_j+A_1$, $2 \leq j \leq s$, and $w=A_1 +B_1$.
Then we have a subgroup
$H=\la\seq{\S}{1}{n}\ra \subset \bar H_2(X;\ZZ)= H_2(X;\ZZ)/\text{torsion}$,
such that $2\bar H_2(X;\ZZ)\subset H$, with $\S_j^2=0$
and $w \cdot \S_j \equiv 1 \pmod 2$, $1 \leq j\leq n$.

We represent $\S_j$ by embedded surfaces of genus $g_j \geq 1$.
Lemma~\ref{lem:1} implies then
\begin{equation}
\label{eqn:rs}
  \DD^w_X(t_1\S_1 +\cdots+ t_n \S_n)=e^{Q(t_1\S_1 +\cdots+ t_n \S_n)/2}
 \vspace{-5mm}  \sum_{-(g_j-1)\leq r_j \leq g_j-1\atop 1\leq j\leq n}
  \vspace{-5mm} a_{r_1\ldots r_n,w} \, e^{2 r_1t_1  +\cdots+2 r_n t_n},
\end{equation}
for complex numbers $a_{r_1\ldots r_n,w}$. They are rational because of the
rationality of the Donaldson invariants. The claim follows.
$K_i$ are integral cohomology classes since $2\bar H_2(X;\ZZ)\subset H$.

\noindent {\underline{Step 3}}.
$K_i$ are lifts of $w_2(X)$ to integer coefficients.

Equivalently, we need to prove that
$K_i \cdot x \equiv x^2 \pmod 2$, for any $x \in \bar H_2(X;\ZZ)$.
Clearly formula~\eqref{eqn:rs} implies that
$K_i \cdot \S_j \equiv 0 = \S_j^2 \pmod 2$, for $1 \leq j\leq n$. Take
$x \in \bar H_2(X;\ZZ) - H$. There is some index $k$ such that $x \cdot \S_k
\neq 0$. We can find an integer $m$ such that
$x' = x+ m\S_k$ has $N=(x')^2 \geq 0$ and $w \cdot x' \equiv 1 \pmod 2$
(recall that $w \cdot \S_j \equiv 1 \pmod 2$ for all $j$).
We blow up $X$ at $N$ points to
get $\tilde X=X \# N \overline{\CP}^2$ with exceptional divisors
$\seq{E}{1}{N}$. Then $y=x' -E_1 - \ldots -E_N$ has $y^2=0$ and $w\cdot y
\equiv 1\pmod 2$.
The blow-up formula~\cite{FS-bl} says
$$
  \DD_{\tilde X}^w= \DD_X^w \cdot e^{-(E_1^2+ \cdots +E_N^2)/2}
  \cosh E_1 \cdots \cosh E_N =
   e^{Q /2} \sum {a_{i,w}\over 2^N}
  \, e^{K_i +\sum \pm E_l},
$$
so the basic classes of $\tilde X$ are of the form $K_i +\sum \pm E_l$, with
$K_i$ basic classes for $X$. Now lemma~\ref{lem:1} applied to $w$ and $y$
implies
in particular that
$$
  0 \equiv (K_i +\sum \pm E_l)\cdot y \equiv K_i \cdot x'  +N \pmod 2
$$
and hence $K_i \cdot x  \equiv x^2 \pmod 2$, for all $K_i$.

\noindent {\underline{Step 4}}.
$\DD^{w'}_X=e^{Q/2} \sum (-1)^{{K_i \cdot {w'} +{w'}^2}\over 2}
a_i e^{K_i}$, for any other $w' \in H^2(X;\ZZ)$, where we put
$a_i= (-1)^{-{{K_i \cdot {w} +{w}^2}\over 2}} a_{i,w}$ (by step 3 the
exponent is an integer).

Lemma~\ref{lem:1} implies that if we have proved the claim for $w'$ and
$w' \cdot \S \equiv 1 \pmod 2$ and $\S^2=0$, then the claim is true
for $w' +\S$. Now to prove the assertion for any $w' \in H^2(X;\ZZ)$,
it is enough to consider $w'=w+L$, with $N=L^2>0$ as in step 1.
Keep those notations and suppose for instance that we are in the case
$(w+E_1)\cdot L' \equiv 1 \pmod 2$
(i.e. $w\cdot L$ is even). Then the blow-up formula again says
$$
  \DD^{w+E_1}_{\tilde X}= -\DD_X^w\cdot e^{-(E_1^2+ \cdots + E_N^2)/2}
  \sinh E_1 \cosh E_2 \cdots \cosh E_N,
$$
so the coefficient of the basic class $K_i +\sum a_l E_l$, where
$(a_l)_{l=1}^N$ is
a sequence of numbers $a_l=\pm 1$, is $(-1)^{a_1+1\over 2} a_{i,w}/2^N$.
By lemma~\ref{lem:1}, 
$$
\DD^{w+L-E_2-\ldots -E_N}_{\tilde X}= e^{Q/2} \sum
c_{K_i +\sum a_l E_l} \left( a_{i,w}/2^N \right) e^{K_i +\sum a_l E_l},
$$ 
where
$$
 c_{K_i +\sum a_l E_l} =(-1)^{a_1+1\over 2} (-1)^{{(K_i +\sum a_l E_l)\cdot
L'\over 2}+1}=
  (-1)^{K_i \cdot L +L^2 \over 2} (-1)^{a_2+1\over 2}\cdots (-1)^{a_N+1\over
2}(-1)^{N-1}.
$$
Using the blow-up formula again together with the standard fact
$D_X^{w+2\a}=(-1)^{\a^2}D_X^w$,
we get $\DD_X^{w+L}=e^{Q/2} \sum (-1)^{K_i \cdot L +2 w\cdot L+L^2 \over 2}
a_{i,w} e^{K_i}$ (since $w\cdot L$ is even), as required.

\noindent {\underline{Step 5}}.
The final assertion is proven as follows. Let $S \inc X$ be an embedded
surface of genus $g$ with $N=S^2 \geq 0$ and representing a non-torsion
homology class. The argument in~\cite[page 709]{KM}
reduces to prove only the case $N>0$. If the genus is $g=0$ then the
Donaldson
invariants vanish identically and hence it is trivially true. In the case
$g \geq 1$, blow-up $X$ at $N$ points to
get $\tilde X=X \# N \overline{\CP}^2$ with exceptional divisors
$\seq{E}{1}{N}$. Consider the proper transform $\tilde S$ of $S$ which
is an embedded surface in $\tilde X$ of genus $g$ 
representing the homology class $S -E_1 - \ldots -E_N$. 
Then lemma~\ref{lem:1} applied to $\tilde X$, $\tilde S$ and $w=E_1$
gives
$$
  2(g-1) \geq | (K_i +\sum \pm E_l)  \cdot (S -E_1 - \ldots -E_N) |
$$
for all $K_i$ basic classes of $X$. Therefore we have $2g-2 \geq
|K_i \cdot S|+S^2$. \hfill $\Box$

\section{Basic classes}
\label{sec:basic}

We first prove theorem~\ref{thm:orderft}. We start with the following
analogue of lemma~\ref{lem:2}.

\begin{lem}
\label{lem:3}
  In the situation of lemma~\ref{lem:1}, $X$ is of
  $w$-finite type of order $k$ if and only if it is of $(w+\S)$-finite type of
order $k$.
\end{lem}

\begin{pf}
  Arguing by contradiction, suppose $X$ is of $w$-finite type of order $k$ and
  of $(w+\S)$-finite type of order $k+a$, $a>0$.
  Then there exists $z=(x^2-4)^{k+a-1} \d_1 \cdots \d_p$ with $\d_i \in
H_1(X)$, $p\geq 0$,
  such that $D^w_X(z e^{tD +s\S+\l x})=0$,
  $D^{w+\S}_X(z e^{tD +s\S+\l x})$ is non-zero but
  $D^{w+\S}_X(z (x^2-4) e^{tD +s\S+\l x})=0$ and
  $D^{w+\S}_X(z \d e^{tD +s\S+\l x})=0$, for any $\d \in H_1(\S)$.
  The arguments in the proof of lemma~\ref{lem:2} now carry over verbatim.
\end{pf}

\noindent {\em Proof of theorem~\ref{thm:orderft}.\/}
It goes as in step 1 of the proof of theorem~\ref{thm:sst} with the difference
that
we use lemma~\ref{lem:3} and note the following:
Let $S$ be a $4$-manifold with $b^+>1$ and $w \in H^2(S;\ZZ)$. Let $\tilde S=
S\# \overline{\CP}^2$ denote its blow-up and let $E$ stand for the cohomology
class of the exceptional divisor. Then $S$ is
of $w$-finite type of order $k$ if and only if $\tilde S$ is of $w$-finite type
of order $k$
if and only if $\tilde S$ is of $(w+E)$-finite type of order $k$. \hfill $\Box$

\noindent {\em Proof of theorem~\ref{thm:basic}.\/}
Now we proceed to the proof of theorem~\ref{thm:basic}. The analogue of
lemma~\ref{lem:1} is the following:

\begin{lem}
\label{lem:4}
  Let $X$ be a $4$-manifold with $b^+>1$ and  $w,\S \in H^2(X;\ZZ)$
  such that $\S^2=0$ and $w \cdot \S \equiv 1 \pmod 2$.
  Let $z\in \L^* H_1(X)$ homogeneous. Take $D\in H_2(X)$. Then
$$
  D^w_X (ze^{tD +s\S+\l x})= e^{Q(tD +s\S)/2+2\l}\sum\limits_{r=-(g-1)}^{g-1}
P_r(t,s,\l) +
  e^{-Q(tD +s\S)/2-2\l}\sum\limits_{r=-(g-1)}^{g-1} Q_r(t,s,\l),
$$
  where $P_r(t,s,\l)$ is a solution of the differential equations $\left(
  {\bd \over \bd s}-2r \right)^N$, $\left( {\bd \over \bd \l} \right)^N$,
  and $Q_r(t,s,\l)$ is a solution of $\left(
  {\bd \over \bd s}-2r\ima \right)^N$, $\left( {\bd \over \bd \l} \right)^N$,
  for $N$ sufficiently large. Moreover
 \begin{eqnarray*}
  D^{w+\S}_X (ze^{tD +s\S+\l x}) &=& e^{Q(tD +s\S)/2+2\l}
  \sum\limits_{r=-(g-1)}^{g-1} (-1)^{r+1} P_r(t,s,\l) + \\
   & & +
   e^{-Q(tD +s\S)/2-2\l}\sum\limits_{r=-(g-1)}^{g-1} (-1)^{r}Q_r(t,s,\l).
  \end{eqnarray*}
\end{lem}

\begin{pf}
  We proceed as in the proof of lemma~\ref{lem:1}, but in this case we have
  $\p^w(X_1, z e^{tD_1+s\S+\l x})\in \whff\subset HFF^*_g$. Now the eigenvalues
  of $(\a, \b,\g)$ in $\whff$ given by proposition~\ref{prop:HFF}
  yield that
$$
  \Dws_X (ze^{tD +s\S+\l x})= e^{Q(tD +s\S)/2+2\l}
  \hspace{-5mm} \sum\limits_{-(g-1) \leq r \leq g-1 \atop
  \text{$r$ odd}} \hspace{-5mm} P_r(t,s,\l) +
  e^{-Q(tD +s\S)/2-2\l}  \hspace{-5mm}
  \sum\limits_{-(g-1) \leq r \leq g-1 \atop
  \text{$r$ even}} \hspace{-5mm} Q_r(t,s,\l),
$$
  where $P_r$ and $Q_r$ satisfy the differential equations in the statement.
Then
  put $d_0=d_0(w) = -w^2 - {3\over 2}(1-b_1+b^+)$, so
  $$
   D^w_X(e^{tD +s\S+\l x})= {1\over 2} \Big(
    e^{Q(tD +s\S)/2+2\l} \big( \sum\limits_{-(g-1) \leq r \leq g-1 \atop
     \text{$r$ odd}}  \hspace{-5mm} P_r(t,s,\l) + \ima^{-d_0} \hspace{-5mm}
   \sum\limits_{-(g-1) \leq r \leq g-1 \atop
  \text{$r$ even}}  \hspace{-5mm} Q_r(\ima t,\ima s,-\l) \big)
$$
$$
 +e^{-Q(tD +s\S)/2-2\l}\big( \sum\limits_{-(g-1) \leq r \leq g-1 \atop
  \text{$r$ even}} \hspace{-5mm} Q_r(t,s,\l) +\ima^{-d_0} \hspace{-5mm}
 \sum\limits_{-(g-1) \leq r \leq g-1 \atop
     \text{$r$ odd}} \hspace{-5mm} P_{-r}(\ima t,\ima s,-\l) \big) \Big).
$$
  We leave $D^{w+\S}_X$ to the reader.
\end{pf}

\noindent {\underline{Step 1}}.
Suppose $X$ is a $4$-manifold with $b^+>1$, $w \in H^2(X;\ZZ)$, $z\in
\L^*H_1(X)$
homogeneous and that there is a subgroup $H=\la\seq{\S}{1}{n}\ra \subset
\bar H_2(X;\ZZ)$
such that $2\bar H_2(X;\ZZ)\subset H$, with $\S_j^2=0$
and $w \cdot \S_j \equiv 1 \pmod 2$, $1 \leq j\leq n$. Then~\eqref{eqn:falta}
holds.

We represent $\S_j$ by embedded surfaces of genus $g_j \geq 1$.
Lemma~\ref{lem:4} implies that
$$
  D^w_X(ze^{t_1\S_1 +\cdots+ t_n \S_n+\l x})=e^{Q(t_1\S_1 +\cdots+ t_n
\S_n)/2+2\l}
 \hspace{-5mm}  \sum_{-(g_j-1)\leq r_j \leq g_j-1\atop 1\leq j\leq n}
  \hspace{-5mm} P_{r_1\ldots r_n,w}(\seq{s}{1}{n},\l) \, e^{2 r_1t_1  +\cdots+2
r_n t_n} +
$$
$$
  e^{-Q(t_1\S_1 +\cdots+ t_n \S_n)/2-2\l}
 \hspace{-5mm}  \sum_{-(g_j-1)\leq r_j \leq g_j-1\atop 1\leq j\leq n}
  \hspace{-5mm} Q_{r_1\ldots r_n,w}(\seq{s}{1}{n},\l) \, e^{\ima (2 r_1t_1
+\cdots+2 r_n t_n)},
$$
where $P$'s and $Q$'s are polynomials. This finishes the step, noting again
that
$K_i$ are integral cohomology classes since $2\bar H_2(X;\ZZ)\subset H$.

\noindent {\underline{Step 2}}.
Fix an arbitrary
$w \in H^2(X;\ZZ)$. Then~\eqref{eqn:falta} holds for $X$ and $w$.

First note that if $w,\S \in H^2(X;\ZZ)$ with
$w\cdot \S\equiv 1\pmod 2$ and $\S^2=0$ then~\eqref{eqn:falta} holds for $w$
implies~\eqref{eqn:falta} holds for $w+\S$.

The relationship of the Donaldson invariants of a $4$-manifold $S$ and the
blow-up $\tilde S=
S\# \overline{\CP}^2$ is not as straightforward as in the simple type case,
therefore the
strategy of proving that~\eqref{eqn:falta} holds for one $w \in H^2(X;\ZZ)$
implies that it holds for any other $w'\in H^2(X;\ZZ)$ by
blowing-up does not work easily. 
Instead we fix $w \in H^2(X;\ZZ)$ and shall look for
a blow-up $\tilde X= X\# m\overline{\CP}^2$ with exceptional divisors
$\seq{E}{1}{m}$ such that there exists $\tilde w \in H^2(\tilde X;\ZZ)$ of the
form
$\tilde w=w +\sum a_i E_i$, $a_i$ integers, with $(\tilde X, \tilde w)$
satisfying the
conditions in step 1. This is enough since $D^w_X(z e^{tD+\l x})=
D^w_{\tilde X}(z e^{tD+\l x})$ and $D^w_X(z e^{tD+\l x})={\bd \over \bd
r}|_{r=0}
D^{w+E}_{\tilde X} (z e^{tD+rE+\l x})$, for $D\in H_2(X)$.

We can blow up to ensure that $Q=r (1) \oplus s(-1)$, $r\geq 2$, $s\geq 1$,
with corresponding basis
$\seq{A}{1}{r},\seq{B}{1}{s}$ and $w \cdot B_n \equiv 0 \pmod 2$.
There are several cases:
\begin{itemize}
  \item Suppose that (up to reordering) $w\cdot A_1 \equiv 1 \pmod 2$,
    $w\cdot A_2 \equiv 0 \pmod 2$. Then blow up once, $\tilde X= X\#
\overline{\CP}^2$,
    with $E$ the exceptional divisor. Set
    \begin{eqnarray*}
     \S_j &=& \left\{ \begin{array}{ll} A_j + E \qquad & \text{if } w\cdot A_j
\equiv 0 \pmod 2
      \\ A_j +B_n &  \text{if }w\cdot A_j \equiv 1 \pmod 2 \end{array} \right.
\\
     \S_{r+j} &=& \left\{ \begin{array}{ll} B_j - A_1 & \text{if } w\cdot B_j
\equiv 0 \pmod 2
      \\ B_j -A_2 &  \text{if }w\cdot B_j \equiv 1 \pmod 2 \end{array} \right.
\\
     \S_{n+1} &=& E-A_2
    \end{eqnarray*}
    Then the subgroup $H=\la\seq{\S}{1}{n+1}\ra \subset H_2(\tilde X;\ZZ)$ and
$\tilde w=w+E$
    satisfy the required properties in step 1.
  \item Suppose that $w\cdot A_j \equiv 1 \pmod 2$, for all $j$. Then put
$\S=A_1+B_n$
    and $w'=w+\S$. Then $w'$ satisfies the conditions of the previous case
so~\eqref{eqn:falta}
   holds for $w'$. Now $w' \cdot \S \equiv 1\pmod 2$ and $\S^2=0$,
so~\eqref{eqn:falta} also holds for $w$.
  \item Suppose that $w\cdot A_j \equiv 0 \pmod 2$, for all $j$ and that there
exists
    $B_1$ such that $w \cdot B_1 \equiv 1\pmod 2$. Then put $\S=A_1+B_1$ and
$w'=w+\S$, and work as in the previous case.
  \item Suppose $w\equiv 0 \pmod 2$. Blow-up once and put $\tilde w=w+E$. It
    reduces to the previous case.
\end{itemize}

Before carrying on with the proof of theorem~\ref{thm:basic}, let us pause to
give a
characterization of the basic classes of $X$. So far
we only can say that the basic classes are relative to a particular $w\in
H^2(X;\ZZ)$. So we define a basic class for $(X,z,w)$ to be a cohomology class
$K_i$ provided by step 2, such that $p_i\neq 0$.

\begin{prop}
\label{prop:charact}
 Let $X$ be a $4$-manifold with $b^+>1$, $z \in \L^*H_1(X)$ homogeneous and
$w\in H^2(X;\ZZ)$.
 Then $K \in H^2(X;\ZZ)$ is a basic class for $(X,z,w)$ if and only if there
exists
 $z'\in \AA_{\text{even}}(X)= \Sym^* (H_0(X)\oplus H_2(X))$ such that
 $D^w_X(z z' e^{tD+\l x})=e^{Q(tD)/2+2\l+K \cdot tD}$.
\end{prop}

\begin{pf}
 Take $N$ larger than the order of finite type of $X$. Then
 $D^w_X(z (1+{x\over 2})^N e^{tD+\l x})=e^{Q(tD)/2+2\l} \sum p_i(tD, \l) e^{K_i
\cdot tD}$.
 Consider a basis $\{ \seq{D}{1}{n}\}$ of $\bar H_2(X;\ZZ)$. Let $\a_j = K
\cdot D_j$, $j=1,
 \ldots, n$. Then take
  $$
  z'= \prod_j \prod_{\b_j \neq \a_j} (D_j -\b_j)^N,
  $$
 for a large enough $N$ (just take it larger than the degrees of the
polynomials $p_i$). Then
 $D^w_X(z (1+{x\over 2})^N z' e^{tD+\l x})=e^{Q(tD)/2+2\l} p(tD, \l) e^{K \cdot
tD}$,
 for a polynomial $p$. Take some $z''=c \,(1-{x\over 2})^m (D_1-\a_1)^{m_1}
\cdots
 (D_n-\a_n)^{m_n}$, for appropriate exponents and constant $c$, to get
 $D^w_X(z (1+{x\over 2})^N z'z'' e^{tD+\l x})=e^{Q(tD)/2+2\l} e^{K \cdot tD}$.
 The converse is obvious.
\end{pf}

We leave the following characterization of basic classes in terms of the
blow-up to the reader.
\begin{cor}
\label{cor:charact}
 Let $X$ be a $4$-manifold with $b^+>1$, $z \in \L^*H_1(X)$ homogeneous and
 $w\in H^2(X;\ZZ)$. Let $\tilde X=X\#\overline{\CP}^2$ be its blow-up with
 exceptional divisor $E$.
 Then $K \in H^2(X;\ZZ)$ is a basic class for $(X,z,w)$ if and only if there
 exists
 $z'\in \AA_{\text{even}}(X)$ such that
 $D^w_{\tilde X}(z z' e^{tD+\l x})=e^{Q(tD)/2+2\l+K \cdot tD} \cosh (E\cdot
tD)$
 if and only if there exists
 $z'\in \AA_{\text{even}}(X)$ such that
 $D^{w+E}_{\tilde X}(z z' e^{tD+\l x})=-e^{Q(tD)/2+2\l+K \cdot tD}
 \sinh (E\cdot tD)$. \hfill $\Box$
\end{cor}

Now we resume the proof of theorem~\ref{thm:basic}.

\noindent {\underline{Step 3}}. The basic classes of $(X,z)$ are independent of
$w\in H^2(X;\ZZ)$.

By lemma~\ref{lem:4} the basic classes for $(X,z,w)$ and $(X,z,w+\S)$ are the
same for any $\S$
such that $\S^2=0$ and $w \cdot \S\equiv 1 \pmod2$. The argument runs as in
step 4 of
section~\ref{sec:sst} using the characterization of basic classes gathered in
corollary~\ref{cor:charact}. We need to use lemma~\ref{lem:4} in the blow-up
manifold
$\tilde X$ with an extra $z' \in \AA_{\text{even}}(\S)$,
but it is easy to see that its statement still holds.

\noindent {\underline{Step 4}}. Any basic class $K_i$ is a lift of $w_2(X)$ to
integer
coefficients.

This is as in step 3 of section~\ref{sec:sst} by noting that
corollary~\ref{cor:charact}
implies in particular that $K \pm E$ are basic classes of the blow-up of $X$.
\hfill
$\Box$

Note also the following corollary to corollary~\ref{cor:charact}

\begin{cor}
 If $K \in H^2(X;\ZZ)$ is a basic class for $X$ then $K \pm E$ are basic
 classes for the blow-up $\tilde X=X \#\overline{\CP}^2$.
\end{cor}

\begin{rem}
  It is natural to expect that there are more basic classes $K +m E$, $m$ odd,
  $m \neq \pm 1$.
\end{rem}

\noindent {\em Acknowledgements:\/} Thanks to Ignasi Mundet for comments on
early versions of this manuscript. I am grateful to the Department of
Mathematics of Universidad Aut\'onoma de Madrid for their hospitality.

\end{document}